\documentclass[3p,12pt,times,fleqn]{elsarticle}
\usepackage[bitstream-charter]{mathdesign}
\usepackage{amsmath}
\usepackage{algorithm}
\usepackage{url}

\def\RR{{\mathbb R}}
\def\CC{{\mathbb C}}
\def\CC{{\mathbb C}}

\def\diag{\mathop{\rm diag}\nolimits}

\newcommand{\ph}{\phantom}

\usepackage{tikz}
\usetikzlibrary{positioning,shapes,arrows}
\usetikzlibrary{trees}

\newtheorem{thm}{Theorem}
\newtheorem{rem}[thm]{Remark}

\newtheorem{lem}[thm]{Lemma}
\newproof{pf}{Proof}
\newtheorem{exm}[thm]{Example}
\newtheorem{alg}[thm]{Algorithm}

\journal{}

\begin{document}
\begin{frontmatter}
\title{Minimal determinantal representations of bivariate polynomials}

\author{Bor~Plestenjak}\ead{bor.plestenjak@fmf.uni-lj.si}
\address{IMFM and Faculty of Mathematics and Physics, University of Ljubljana, Jadranska 19,
  SI-1000 Ljubljana, Slovenia.}

\begin{abstract} For a square-free bivariate polynomial $p$ of degree $n$ we
introduce a simple and fast numerical algorithm
for the construction of $n\times n$ matrices $A$, $B$, and $C$ such that
$\det(A+xB+yC)=p(x,y)$. This is the minimal size needed
to represent a bivariate polynomial of degree $n$. Combined with a
square-free factorization one can now
compute $n \times n$ matrices for any bivariate polynomial of degree $n$.
The existence of such symmetric matrices was established by Dixon in 1902, but, up to
now, no simple numerical construction has been found, even if
the matrices can be  nonsymmetric.
Such representations may be used to efficiently numerically solve a system
of two bivariate polynomials of small degree
via the eigenvalues of a two-parameter eigenvalue problem.
The new representation speeds up the computation considerably.\\

\end{abstract}

\end{frontmatter}

\tikzstyle{tocka} = [circle, draw, fill=blue!40, inner sep=0pt, minimum size=6pt, node distance=1em]
\tikzstyle{tockar} = [circle, draw, fill=red!80, inner sep=0pt, minimum size=6pt, node distance=1em]
\tikzstyle{tockaw} = [circle, draw, fill=white!80, inner sep=0pt, minimum size=6pt, node distance=1em]

\section{Introduction}
Let
\begin{equation}
p(x,y) := \sum_{i=0}^n \, \sum_{j=0}^{n-i} \ p_{ij} \, x^i \, y^j,\label{eq:pol}
\end{equation} where
$p_{ij}\in\CC$ for all $i,j$, be a bivariate polynomial of degree $n$,
where we assume that $p_{ij}\ne 0$ for at least one index such that $i+j=n$.
We say that
matrices $A,B,C\in\CC^{m\times m}$, where $m\ge n$, form a
\emph{determinantal representation} of order $m$ of the
polynomial $p$ if
\begin{equation}\label{eq:detr}
\det(A+xB+yC)=p(x,y).
\end{equation}

It is known since Dixon's 1902 paper \cite{Dixon} that every
bivariate polynomial of degree $n$ admits
a determinantal representation with symmetric matrices of
order $n$.
However, the construction of such matrices is far from trivial and
up to now there have been no efficient numerical algorithms, even if
we do not insist on matrices being symmetric.
We introduce the first efficient numerical construction
of determinantal representations that returns
$n\times n$ matrices for a square-free bivariate polynomial of degree $n$,
which, with the exception of the symmetry, agrees with Dixon's result.
For non square-free
polynomials  one can combine it with a
square-free factorization to obtain a representation of order $n$.

Our motivation comes from the following approach for finding roots of systems of
bivariate polynomials, proposed by Plestenjak and Hochstenbach in \cite{BorMichiel}.
Suppose that we have a system of two bivariate polynomials
\begin{equation}\label{eq:p}
\begin{split}
p(x,y) &:= \sum_{i=0}^{n_1} \, \sum_{j=0}^{{n_1}-i} \ p_{ij} \, x^i \, y^j=0, \\
q(x,y) &:= \sum_{i=0}^{n_2} \, \sum_{j=0}^{{n_2}-i} \ q_{ij} \, x^i \, y^j=0.
\end{split}
\end{equation}
The idea is to construct matrices
$A_1,B_1,C_1$ of size $m_1\times m_1$ and matrices $A_2, B_2, C_2$ of
size $m_2\times m_2$ such that
\begin{equation}\label{eq:dr}
\begin{split}
\det(A_1+x B_1 +y C_1)&=p(x,y),\\[1mm]
\det(A_2+x B_2 +y C_2)&=q(x,y)
\end{split}
\end{equation}
and then numerically solve the equivalent two-parameter eigenvalue problem \cite{Atkinson}
\begin{equation}\label{eq:par2}
\begin{split}
(A_1+x B_1 +y C_1) \, u_1&=0,\\[1mm]
(A_2+x B_2 +y C_2) \, u_2&=0.
\end{split}
\end{equation}
Here it does not matter whether the matrices are symmetric. Except in some special cases,
we get $m_1m_2> n_1n_2$ and thus \eqref{eq:par2}
is a singular two-parameter eigenvalue problem.
Its finite regular eigenvalues, which can be computed numerically with a staircase type
algorithm \cite{MuhicPlestenjakLAA}, are the roots of \eqref{eq:p}. This is a generalization
of a well-known approach to numerically compute the roots of a univariate polynomial as
eigenvalues of its companion matrix.

The first phase, where we construct determinantal representations, should be
computationally simple and return matrices as small as possible. It is
easy to see that for a polynomial of degree $n$ we need matrices
of order at least $n$.
There exist procedures, which use slow symbolic computation and other expensive steps, that do
return representations of order $n$ for certain bivariate polynomials,
but they are computationally too complex. For example, the method from \cite{Plaumann},
based on the proof from \cite{Dixon}, returns $n\times n$ symmetric matrices for
a polynomial of degree $n$ that satisfies the so called real zero condition,
but it is too slow for our purpose. As the first step in \cite{Plaumann}
is to find zeros of a system of bivariate polynomials
of degree $n$ and $n-1$, this clearly can not be efficient enough for our needs.
In addition, we need a determinantal representation for any
bivariate polynomial of degree $n$.

In our approach to the computation of roots of a system of bivariate polynomials we must find the
right balance. As we can not exploit the symmetry
in \eqref{eq:dr}, it is not necessary that the matrices have this property. If we want
a representation with small matrices this might take too much time in the first phase, while a
representation with large matrices slows down the second phase.
By using the most efficient  representations at the time with
 matrices of order $\frac{1}{6}n^2+{\cal O}(n)$, a numerical method
was implemented in Matlab \cite{BorBR} that is competitive
to some existing methods for polynomials up to degree 9, see \cite{BorMichiel} for details.

The computational complexity of the new construction of determinantal representations
that we give in the paper is just slightly above that of the construction used in
\cite{BorMichiel}, yet the matrices are much smaller. While the matrices in \cite{BorMichiel}
are of order $\frac{1}{6}n^2+{\cal O}(n)$, the new construction gives matrices of the minimal
possible order $n$ for a square-free bivariate polynomial of degree $n$. This decreases the
overall asymptotic complexity of solving the system \eqref{eq:p} via determinantal
representations from ${\cal O}(n^{12})$
to ${\cal O}(n^6)$ when $n=n_1=n_2$. Moreover, if both polynomials $p_1$ and $p_2$
are square-free, then
$m_1m_2=n_1n_2$ and \eqref{eq:par2} is a nonsingular two-parameter eigenvalue problem, which is
much easier to solve numerically than a singular one.

The rest of the paper is organized as follows. In Section~\ref{sec:over} we
give a short overview of existing determinantal representations.
In Section~\ref{sec:preproc} we show that a square-free
polynomial can always be transformed
into a form required by the algorithm in Section~\ref{sec:two}, where we
give a determinantal representation of order $n$ for a square-free bivariate polynomial
of degree $n$. In Section~\ref{sec:three} we extend the representation of order $n$ to
non square-free polynomials and discuss other options for such polynomials.
We end with some numerical experiments in Section~\ref{sec:numex}
and conclusions in Section \ref{sec:conc}.

\section{Overview of existing determinantal representations}\label{sec:over}

In the semidefinite programming (SDP) there is a large interest in
symmetric determinantal representations
of the real zero polynomials, a particular subset of polynomials
related to the linear matrix inequality (LMI)
constraints. For an overview, see, e.g., \cite{Netzer, VinnikovSurvey}.
A real bivariate polynomial $p$ satisfies the \emph{real zero condition}
with respect to $(x_0,y_0)\in\RR^2$
if for all $(x,y)\in\RR^2$ the univariate polynomial
$p_{(x,y)}(t):=p(x_0+tx,y_0+ty)$ has only real zeros.
A two-dimensional \emph{LMI set} is defined as
$\left\{(x,y)\in\RR^2: A+ xB+y C\succeq 0\right\}$,
where $A$, $B$, and $C$ are symmetric matrices of size $m\times m$
and $\succeq 0$ stands
for positive semidefinite.
For this particular subset of polynomials there do exist some procedures that involve slow
symbolic computation and other expensive steps,
see, e.g., \cite{Plaumann, PlaumannHVcurves}. However, besides being too slow these algorithms
are limited to the real zero polynomials only.

While in SDP and LMI the matrices
have to be symmetric or Hermitian, this is not important in our case.
We are looking for a simple and fast numerical construction
of matrices as small as possible that satisfy (\ref{eq:detr}).

Here is a list of some available determinantal representations for generic
bivariate polynomials.
The first group of determinantal representations has the property that
the elements of matrices $A$, $B$, and $C$ depend affine-linearly on
the coefficients of the polynomial $p$. Such determinantal representations
are named \emph{uniform} in \cite{Jan}.
The first such representation of order $n^2$ is given by
Khazanov in \cite{Khazanov} as a special case of a linearization of
a multiparameter polynomial matrix.
This is improved to a representation or order $\frac{1}{2} n(n+1)$ by
Muhi\v{c} and Plestenjak in
\cite[Appendix]{MuhicPlestenjakLAA}. Quarez \cite{Quarez} gives
symmetric representations of multivariate polynomials, which results in a representation of order $\frac{1}{4}n^2+{\cal O}(n)$ for a bivariate polynomial.
A smaller nonsymmetric uniform representation of the same asymptotic order is
described in \cite{BorMichiel}. Recently, a  uniform
representation of order $2n-1$ was presented in \cite{Jan}, which
is the first uniform representation such that the order of matrices
grows linearly and not quadratically with $n$.
All uniform representations do not require any computation, the construction
is very simple and fast as one just puts
the coefficients of $p$ on prescribed places in the matrices $A$, $B$, and $C$.

If we allow computations, we can obtain smaller representations.
In \cite{BorMichiel}, a representation
of order $\frac{1}{6}n^2+{\cal O}(n)$  is given.
This representation is used in \cite{BorBR} as a part of a numerical method
for the roots of a system of bivariate polynomials that is competitive to
the existing numerical methods for polynomials of small degree.
In
the following  sections we upgrade the approach from \cite{BorMichiel} to
obtain minimal determinantal representations of order $n$  while
maintaining approximately the same complexity of computations involved in the
construction.

\section{Preliminary transformations}\label{sec:preproc}

We can homogenize \eqref{eq:pol} into
\begin{equation}\label{eq:homop}
p_h(x,y,z):=z^np\left(\frac{x}{z},\frac{y}{z}\right)=
\sum_{i=0}^{n} \, \sum_{j=0}^{{n}-i} \ p_{ij} \, x^i \, y^j z^{n-i-j},
\end{equation}
where $(x,y,z)$ are points in the projective plane.
It is easy to see that if $n\times n$
matrices $A,B,C$ are such that
$\det(A+xB+yC)=p(x,y)$, then
\begin{equation}\label{eq:detreph}
\det(zA+xB+yC)={p_h}(x,y,z).
\end{equation} This also works in the
opposite direction. If we construct a
determinantal representation \eqref{eq:detreph} with matrices
of order $n$
of the homogeneous polynomial $p_h$, then
we get a determinantal representation of $p(x,y)$
by simply setting $z=1$.

The homogeneous form gives us are more freedom in the following sense.
We can apply a linear change of variables
\begin{equation}\label{eq:trt}\left[\begin{matrix}x \cr y \cr z\end{matrix}\right]= T
\left[\begin{matrix}\widetilde{x} \cr \widetilde{y} \cr \widetilde{z}\end{matrix}\right],
\end{equation}
where $T$ is a nonsingular $3\times 3$ matrix, and transform
$p_h(x,y,z)$ into a homogeneous polynomial $\widetilde{p}_h(\widetilde{x},\widetilde{y},\widetilde{z})$. If we find matrices $\widetilde{A}$, $\widetilde{B}$, and $\widetilde{C}$ for $\widetilde{p}_h$ such that
$\det(\widetilde{z}\widetilde{A}+\widetilde{x}\widetilde{B}+\widetilde{y}\widetilde{C})=
{\widetilde{p}_h}(\widetilde{x},\widetilde{y},\widetilde{z})$,
this gives a determinantal representation of $p$ after we substitute
$\widetilde{x},\widetilde{y},\widetilde{z}$ back
to $x,y,z$ and set $z=1$.

We say that a polynomial \eqref{eq:pol} is square-free if it is not a multiple
of a square of a non-constant polynomial. The following result that
 relies on B\'ezout's theorem and Bertini's theorem (see, e.g., \cite{Smith}),
 ensures that our construction
 can be applied to any square-free polynomial.

\begin{thm}\label{thm:fact}
  Let $p$ be a bivariate non-constant polynomial of degree $n$. A
  generic line ${\cal L}$ of the form $\alpha x +\beta y +\gamma z=0$ in the projective plane intersects
  the curve $p_h(x,y,z)=0$, where $p_h$ is the homogenized polynomial $p$, in $n$ distinct
  points if and only if $p$ is square-free.
\end{thm}

\begin{pf}
  Let ${\cal C}$ be the zero set of $p_h$.
  When ${\cal C}$ is smooth, it follows from Bertini's theorem that a generic line ${\cal L}$ intersects ${\cal C}$ in distinct points. This can be easily extended to the case when ${\cal C}$ has finitely many singular points, which is true for a square-free polynomial. Namely, in the Zariski topology the lines containing a singular point form
  a proper closed subset of all possible lines.
  As $p$ has degree $n$, there are $n$ points in ${\cal L}\cap{\cal C}$ by B\'ezout's theorem.

On the other hand, if $p$ is not square-free then $p(x,y)=q(x,y)^kr(x,y)$, where $k\ge 2$ and
$q$ is a non-constant polynomial. Clearly, each intersection of an arbitrary line ${\cal L}$ with the zero set
of $q_h$, the homogenization of $q$, appears with multiplicity at least $k$ in ${\cal L}\cap {\cal C}$. \qed
\end{pf}

Let \eqref{eq:pol} be a bivariate polynomial $p$ of degree $n$ that we want to linearize.
We can assume that $p$
has the following properties:
\begin{enumerate}
\item[(a)] $p_{n0}\ne 0$,
\item[(b)] all zeros $\xi_1,\ldots,\xi_{n}$ of the polynomial
\begin{equation}\label{eq:rootsxi}
h(\xi):={p_h}(\xi,1,0)=p_{n0}\xi^{n}+p_{n-1,1}\xi^{n-1}+\cdots+p_{0n}=0
\end{equation}
are simple.
\end{enumerate}
The above holds for a generic square-free polynomial $p$. If not,
one can apply a random linear substitution \eqref{eq:trt}. In particular,
an equivalent formulation of (b) is that
in the projective plane the line $z=0$ and the curve $p_h(x,y,z)=0$ intersect
in $n$ distinct points. If follows from Theorem \ref{thm:fact} that this holds
after a random linear substitution \eqref{eq:trt}.
\medskip

In a preliminary step of the construction we apply a linear substitution
of the form
\[x=\widetilde x + s\widetilde y + t \widetilde z,\quad y=\widetilde y,
\quad z=\widetilde z,\]
where we set $s$ and $t$ so that
$\widetilde p_{0n}=\widetilde p_{0,n-1}=0$ in
the transformed
polynomial $\widetilde p$ while properties
 (a) and (b) still hold.
Indeed,
$\widetilde p_{n0}=\widetilde p_h(1,0,0)=p_h(1,0,0)=p_{n0}\ne 0$
and it follows from
$\widetilde p_h(\xi,1,0) = p_h(\xi +s, 1, 0)$ that all roots of \eqref{eq:rootsxi} shift for $s$ and thus remain simple.
It is easy to see that
\[\widetilde p_{0n}= \widetilde p_h(0,1,0) =  p_h(s,1,0)=h(s)\]
and
\[\widetilde p_{0,n-1} = {d\over d\widetilde z}\widetilde p_h(0,1,0) =
t {d\over dx}p_h(s,1,0) + {d\over dz}p_h(s,1,0) = t h'(s) + {d\over dz}p_h(s,1,0).\]
Therefore, if we select $s$ as one of the roots of \eqref{eq:rootsxi}
and
\[ t=-\frac{p_{n-1,0}s^{n-1}+p_{n-2,1}s^{n-2}+\cdots + p_{1,n-2}s+p_{0,n-1}}{ h'(s)},\]
then coefficients $\widetilde p_{0n}$ and $\widetilde p_{0,n-1}$ are both zero.
Note that $t$ is well defined because all roots of \eqref{eq:rootsxi} are simple
 and thus $h'(s)\ne 0$.

\section{Determinantal representation for a square-free polynomial}\label{sec:two}

Let \eqref{eq:pol} be a bivariate square-free polynomial of degree $n$ that we want to linearize. After the preliminary linear transformations from the previous section we can assume that:
\begin{enumerate}
\item[(a)] $p_{n0}\ne 0$,
\item[(b)] all zeros $\xi_1,\ldots,\xi_{n-1}$ of the polynomial
\begin{equation}\label{eq:racxi}
v(\xi):=p_{n0}\xi^{n-1}+p_{n-1,1}\xi^{n-2}+\cdots+p_{1,n-1}
\end{equation}
are simple and nonzero,
\item[(c)] $p_{0n}=p_{0,n-1}=0$.
\end{enumerate}

The construction is based on  bivariate polynomials $q_0,\ldots,q_{n-1},$
constructed recursively as
\begin{align}
q_0(x,y) & :=1,\nonumber\\
q_1(x,y) & :=f_{11} q_0(x,y),\nonumber\\
q_2(x,y) & :=f_{21} q_1(x,y) + f_{22} q_0(x,y),\nonumber\\
& \vdots\label{eq:q}\\
q_{n-1}(x,y) & := f_{n-1,1} q_{n-2}(x,y) + f_{n-1,2} q_{n-3}(x,y) + \cdots + f_{n-1,n-1} q_0(x,y),\nonumber
\end{align}
where each coefficient $f_{ij}$ has a linear
form $f_{ij}=\alpha_{ij}x+\beta_{ij}y$ with $\alpha_{ij}=1$ for $i<n-1$. It follows
from the construction \eqref{eq:q} that $q_j$ is a bivariate polynomial of degree $j$ for $j=0,\ldots,n-1$.
The ansatz for a determinantal representation of $p$ is an $n\times n$
bivariate pencil
\begin{equation}\label{eq:ansatz}
A+xB+yC=\left[\begin{matrix}\gamma_{00} + \gamma_{10}x & \gamma_1 & \gamma_2 & \cdots & \gamma_{n-2} & p_{n0}x \cr
-f_{11} & 1 \cr
-f_{22} & -f_{21} & 1\cr
-f_{33} & -f_{32} & - f_{31} & 1\cr
\vdots & & & \ddots & \ddots \cr
-f_{n-1,n-1} & -f_{n-1,n-2} & \cdots & \cdots & -f_{n-1,1} & 1
\end{matrix}\right].
\end{equation}
One can see from \eqref{eq:q} and \eqref{eq:ansatz} that
\[(A+xB+yC)\left[\begin{matrix}1 &  & &  &  \cr
q_1(x,y) & 1 \cr
q_2(x,y) &  & 1\cr
\vdots & &  & \ddots  \cr
q_{n-1}(x,y) &  &  &  & 1\end{matrix}\right]
=
\left[\begin{matrix}d(x,y) & \gamma_1 & \gamma_2 & \cdots & \gamma_{n-2} & p_{n0}x \cr
0 & 1 \cr
0 & -f_{21} & 1\cr
0 & -f_{32} & - f_{31} & 1\cr
\vdots & \vdots & & \ddots & \ddots \cr
0 & -f_{n-1,n-2} & \cdots & \cdots & -f_{n-1,1} & 1
\end{matrix}\right],
\]
where
\begin{equation}\label{eq:expand}
d(x,y)=\gamma_{00} + \gamma_{10}x+\gamma_{1}q_1(x,y) + \cdots
  +\gamma_{n-2} q_{n-2}(x,y) + p_{n0} x q_{n-1}(x,y).
\end{equation}
It follows that $\det(A+xB+yC)=d(x,y)$ and
we will show how to set the values of $\alpha_{ij},\beta_{ij}$,
and $\gamma_{00},\gamma_{10},\gamma_1,\ldots,\gamma_{n-2}$ so
that $\det(A+xB+yC)=p(x,y)$.
\medskip

For a better understanding we first give a quick
overview of the algorithm. Then we explain the details
and give the complete procedure in Algorithm \ref{alg:osijek}, followed by an  example. In the quick overview
we display the structure of a polynomial by a
 diagram, where dots in the $j$-th row stand
for zero or nonzero coefficients at $x^{j-1}, x^{j-2}y, \ldots, y^{j-1}$, respectively.
The following diagram is for the case $n=5$.
Notice the two white dots representing $p_{04}=p_{05}=0$.
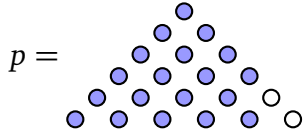
\[
p=
\begin{tikzpicture}[auto,thick,main node/.style={tocka},baseline=(current bounding box.center)]
  \node[tocka] (00) {};
  \node[tocka] (10) [below left of=00] {};
  \node[tocka] (01) [below right of=00] {};
  \node[tocka] (20) [below left of=10] {};
  \node[tocka] (11) [below right of=10] {};
  \node[tocka] (02) [below right of=01] {};
  \node[tocka] (30) [below left of=20] {};
  \node[tocka] (21) [below right of=20] {};
  \node[tocka] (12) [below right of=11] {};
  \node[tocka] (03) [below right of=02] {};
  \node[tocka] (40) [below left of=30] {};
  \node[tocka] (31) [below right of=30] {};
  \node[tocka] (22) [below right of=21] {};
  \node[tocka] (13) [below right of=12] {};
  \node[tockaw] (04) [below right of=03] {};
  \node[tocka] (50) [below left of=40] {};
  \node[tocka] (41) [below right of=40] {};
  \node[tocka] (32) [below right of=31] {};
  \node[tocka] (23) [below right of=22] {};
  \node[tocka] (14) [below right of=13] {};
  \node[tockaw] (05) [below right of=04] {};
\end{tikzpicture}
\]

We build the representation
\eqref{eq:ansatz}
in a loop by adding subdiagonals from top to bottom and
elements in the first row from right to left. We start
by
taking $\alpha_{i1}=1$ and $\beta_{i1}=-\xi_{i}$, i.e., $f_{i1}=x-\xi_{i}y$, for $i=1,\ldots,n-1$. Because of \eqref{eq:racxi} $p_{n0}x (x-\xi_1 y)\cdots (x-\xi_{n-1}y)$ agrees with the part of $p$ of degree $n$ and we get a residual of degree $n-1$ with
a zero coefficient at $y^{n-1}$. In case $n=5$ we get
\begin{equation}\label{eq:slika1}
r^{(1)}(x,y):=p(x,y)-\det\left[\begin{matrix} 0 &  &  &   & p_{50}x \cr
-f_{11} & 1 \cr
 & -f_{21} & 1\cr
&  & -f_{31} & 1\cr
 &  &  &  -f_{41} & 1\end{matrix}\right] =
\begin{gathered}
 \begin{tikzpicture}[auto,node distance=1em,thick,main node/.style={tocka},baseline=(current bounding box.center)]
  \node[tocka] (00) {};
  \node[tocka] (10) [below left of=00] {};
  \node[tocka] (01) [below right of=00] {};
  \node[tocka] (20) [below left of=10] {};
  \node[tocka] (11) [below right of=10] {};
  \node[tocka] (02) [below right of=01] {};
  \node[tocka] (30) [below left of=20] {};
  \node[tocka] (21) [below right of=20] {};
  \node[tocka] (12) [below right of=11] {};
  \node[tocka] (03) [below right of=02] {};
  \node[tocka] (40) [below left of=30] {};
  \node[tocka] (31) [below right of=30] {};
  \node[tocka] (22) [below right of=21] {};
  \node[tocka] (13) [below right of=12] {};
  \node[tockaw] (04) [below right of=03] {};
  \node[tockaw] (50) [below left of=40] {};
  \node[tockaw] (41) [below right of=40] {};
  \node[tockaw] (32) [below right of=31] {};
  \node[tockaw] (23) [below right of=22] {};
  \node[tockaw] (14) [below right of=13] {};
  \node[tockaw] (05) [below right of=04] {};
\end{tikzpicture}
\end{gathered},
\end{equation}
where $r^{(1)}$ is a polynomial of  degree $4$ with a zero coefficient at $y^4$.
Now we add $f_{22},\ldots,f_{n-1,2}$ in the second subdiagonal to
annihilate the part of degree $n-1$. In case $n=5$ we get
\[s^{(2)}(x,y):=p(x,y)-\det\left[\begin{matrix} 0 &  &  &   & p_{50}x \cr
-f_{11} & 1 \cr
 -f_{22} & -f_{21} & 1\cr
&  -f_{32} & -f_{31} & 1\cr
 &  & -f_{42} &  -f_{41} & 1\end{matrix}\right] =
\begin{gathered}
 \begin{tikzpicture}[auto,node distance=1em,thick,main node/.style={tocka},baseline=(current bounding box.center)]
  \node[tocka] (00) {};
  \node[tocka] (10) [below left of=00] {};
  \node[tocka] (01) [below right of=00] {};
  \node[tocka] (20) [below left of=10] {};
  \node[tocka] (11) [below right of=10] {};
  \node[tocka] (02) [below right of=01] {};
  \node[tocka] (30) [below left of=20] {};
  \node[tocka] (21) [below right of=20] {};
  \node[tocka] (12) [below right of=11] {};
  \node[tocka] (03) [below right of=02] {};
  \node[tockaw] (40) [below left of=30] {};
  \node[tockaw] (31) [below right of=30] {};
  \node[tockaw] (22) [below right of=21] {};
  \node[tockaw] (13) [below right of=12] {};
  \node[tockaw] (04) [below right of=03] {};
  \node[tockaw] (50) [below left of=40] {};
  \node[tockaw] (41) [below right of=40] {};
  \node[tockaw] (32) [below right of=31] {};
  \node[tockaw] (23) [below right of=22] {};
  \node[tockaw] (14) [below right of=13] {};
  \node[tockaw] (05) [below right of=04] {};
\end{tikzpicture}
\end{gathered}.
\]
There are $n-1$ free parameters $\beta_{22},\ldots,\beta_{n-1,2}$ and $\alpha_{n-1,2}$ in $f_{22},\ldots,f_{n-1,2}$ to zero
$n-1$ coefficients in the residual.  We show later that this can be done by applying
a direct formula for $\alpha_{n-1,2}$ and solving
a nonsingular triangular system of linear equations for $\beta_{22},\ldots,\beta_{n-1,2}$. Notice that the introduction of $f_{22},\ldots,f_{n-1,2}$ does not affect the part of the determinant of  degree $n$. Next we add $\gamma_{n-2}$ to make
the coefficient at $y^{n-2}$ zero. In case $n=5$
we get
\begin{equation}\label{eq:slika2}
r^{(2)}(x,y):=p(x,y)-\det\left[\begin{matrix} 0 &  &  & \gamma_3  & p_{50}x \cr
-f_{11} & 1 \cr
 -f_{22} & -f_{21} & 1\cr
&  -f_{32} & -f_{31} & 1\cr
 &  & -f_{43} &  -f_{41} & 1\end{matrix}\right] =
\begin{gathered}
 \begin{tikzpicture}[auto,node distance=1em,thick,main node/.style={tocka},baseline=(current bounding box.center)]
  \node[tocka] (00) {};
  \node[tocka] (10) [below left of=00] {};
  \node[tocka] (01) [below right of=00] {};
  \node[tocka] (20) [below left of=10] {};
  \node[tocka] (11) [below right of=10] {};
  \node[tocka] (02) [below right of=01] {};
  \node[tocka] (30) [below left of=20] {};
  \node[tocka] (21) [below right of=20] {};
  \node[tocka] (12) [below right of=11] {};
  \node[tockaw] (03) [below right of=02] {};
  \node[tockaw] (40) [below left of=30] {};
  \node[tockaw] (31) [below right of=30] {};
  \node[tockaw] (22) [below right of=21] {};
  \node[tockaw] (13) [below right of=12] {};
  \node[tockaw] (04) [below right of=03] {};
  \node[tockaw] (50) [below left of=40] {};
  \node[tockaw] (41) [below right of=40] {};
  \node[tockaw] (32) [below right of=31] {};
  \node[tockaw] (23) [below right of=22] {};
  \node[tockaw] (14) [below right of=13] {};
  \node[tockaw] (05) [below right of=04] {};
\end{tikzpicture}
\end{gathered}.
\end{equation}

By comparing \eqref{eq:slika2} to \eqref{eq:slika1} we see that
we reduced the degree of the residual by one and preserved the form where
the coefficient at $y^d$, where $d$ is the degree of the residual, is zero.
In a similar way we fix the remaining parameters. We
start the step with the residual $r^{(k-1)}$, which has degree $n-k+1$ and a zero
coefficient at $y^{n-k+1}$. In the first part we set
$f_{kk},\ldots,f_{n-1,k}$ so that
the new residual $s^{(k)}$ has degree $n-k$. Then
we set $\gamma_{n-k}$ to annihilate the term $y^{n-k}$ in
$s^{(k)}$ and obtain the residual $r^{(k)}$ for the next step.
We continue in the same manner until the final
residual $r^{(n-1)}$ has the form $\gamma_{00} + \gamma_{10} x$. In our example, after two additional steps, we get
\[
r^{(4)}(x,y):=p(x,y)-\det\left[\begin{matrix} 0 & \gamma_1 & \gamma_2 & \gamma_3  & p_{50}x \cr
-f_{11} & 1 \cr
 -f_{22} & -f_{21} & 1\cr
  -f_{33} & -f_{32} & -f_{31} & 1\cr
 -f_{44} & -f_{43} & -f_{42} &  -f_{41} & 1\end{matrix}\right] =
\begin{gathered}
 \begin{tikzpicture}[auto,node distance=1em,thick,main node/.style={tocka},baseline=(current bounding box.center)]
  \node[tocka] (00) {};
  \node[tocka] (10) [below left of=00] {};
  \node[tockaw] (01) [below right of=00] {};
  \node[tockaw] (20) [below left of=10] {};
  \node[tockaw] (11) [below right of=10] {};
  \node[tockaw] (02) [below right of=01] {};
  \node[tockaw] (30) [below left of=20] {};
  \node[tockaw] (21) [below right of=20] {};
  \node[tockaw] (12) [below right of=11] {};
  \node[tockaw] (03) [below right of=02] {};
  \node[tockaw] (40) [below left of=30] {};
  \node[tockaw] (31) [below right of=30] {};
  \node[tockaw] (22) [below right of=21] {};
  \node[tockaw] (13) [below right of=12] {};
  \node[tockaw] (04) [below right of=03] {};
  \node[tockaw] (50) [below left of=40] {};
  \node[tockaw] (41) [below right of=40] {};
  \node[tockaw] (32) [below right of=31] {};
  \node[tockaw] (23) [below right of=22] {};
  \node[tockaw] (14) [below right of=13] {};
  \node[tockaw] (05) [below right of=04] {};
\end{tikzpicture}
\end{gathered},
\]
where $r^{(4)}(x,y)=\gamma_{00} + \gamma_{10}x$.
We obtain the final determinantal representation by putting
$\gamma_{00} + \gamma_{10}x$ at position $(1,1)$ in \eqref{eq:ansatz}.
\medskip

Let us show in more details the step from the residual $r^{(k-1)}$ to the residual $r^{(k)}$. Let $q_j^{(k-1)}$ for $j=0,\ldots,n-1$ be the polynomials \eqref{eq:q} in step $k-1$, i.e., $q_0^{(k-1)}(x,y)=1$ and
\begin{equation}\label{eq:qdif}
q_j^{(k-1)}(x,y)=\sum_{\ell=1}^{\min(j,k-1)} f_{j\ell}q_{j-\ell}^{(k-1)}(x,y),\qquad j=1,\ldots,n-1.\end{equation}
It is clear from \eqref{eq:qdif} that
$q_j^{(k-1)}$ and $q_j^{(k)}$ are equal for $j<k$ and
differ only in terms of degree $j-k+1$ or less for $j\ge k$.
We take
\begin{align}
r^{(k-1)}(x,y)&:=p(x,y)-p_{n0}x q_{n-1}^{(k-1)}(x,y)-\gamma_{n-2}q_{n-2}^{(k-1)}(x,y)-\cdots-
 \gamma_{n-k+1}q_{n-k+1}^{(k-1)}(x,y)\label{eq:respk1}\\
&\phantom{:}=r^{(k-1)}_{00}+r^{(k-1)}_{10}x+r^{(k-1)}_{01}y+\cdots+r^{(k-1)}_{n-k+1,0}x^{n-k+1}+
  \cdots+r^{(k-1)}_{1,n-k}x y^{n-k},\nonumber
\end{align}
select the terms of
 degree $n-k+1$ and
form
\[u_{n-k+1}(x,y):=r^{(k-1)}_{n-k+1,0}x^{n-k+1}+ \cdots+r^{(k-1)}_{1,n-k}x y^{n-k}.\]
Notice that $r_{0,n-k+1}^{(k-1)}=0$, therefore we can write
\begin{equation}\label{eq:polh}
u_{n-k+1}(x,y) = p_{n0} x h_{n-k}(x,y),
\end{equation}
where $h_{n-k}$ is a polynomial of degree $n-k$.
It follows
  from \eqref{eq:respk1} and \eqref{eq:polh} that in order to zero all terms of degree $n-k+1$ in $r^{(k-1)}$ the part of degree $n-k$ in $q_{n-1}^{(k)}-q_{n-1}^{(k-1)}$ has
to agree with $h_{n-k}$.
By comparing
$q_{n-1}^{(k)}$ to $q_{n-1}^{(k)}$ we see that we have to set
the parameters $\beta_{kk},\ldots,\beta_{n-1,k}$ and $\alpha_{n-1,k}$ so
that
\begin{align*}
h_{n-k}(x,y)&=\sum_{\ell=k}^{n-1}f_{\ell k}\prod_{i=1}^{\ell-k}f_{i1} \prod_{j=\ell+1}^{n-1}f_{j1} \\
&=f_{k,k}f_{k+1,1}\ldots f_{n-1,1} + f_{11}f_{k+1,k}f_{k+2,1}\ldots f_{n-1,1} + \cdots +
f_{11}\ldots f_{n-k-1,1}f_{n-1 k}.
\end{align*}
This is equivalent to finding $\beta_{kk},\ldots,\beta_{n-1,k}$, and $\alpha_{n-1,k}$ such
that
\begin{align}
h_{n-k}(t,1)&=\sum_{\ell=k}^{n-2}(t+\beta_{\ell k})\prod_{i=1}^{\ell-k}(t-\xi_{i}) \prod_{j=\ell+1}^{n-1}(t-\xi_{j})
+(\alpha_{n-1, k}t+\beta_{n-1,  k})\prod_{i=1}^{n-k-1}(t-\xi_{i}).\label{eq:razh}
\end{align}
By inspecting the coefficients at $t^{n-k}$ in \eqref{eq:razh} we get
\begin{equation}\label{eq:alpha}
\alpha_{n-1,k}=h_{n-k}(1,0)-n+k+1.
\end{equation}
For the remaining parameters $\beta_{kk},\ldots,\beta_{n-1,k}$ it follows from \eqref{eq:razh}
that
\begin{align}
g_{n-k-1}(t)&:=h_{n-k}(t,1)- t
\sum_{\ell=k}^{n-2}\prod_{i=1}^{\ell-k}(t-\xi_{i}) \prod_{j=\ell+1}^{n-1}(t-\xi_{j})
-\alpha_{n-1, k} t\prod_{i=1}^{n-k-1}(t-\xi_{i})\label{eq:eng1}\\
&\phantom{:}=\sum_{\ell=k}^{n-2}\beta_{\ell k}\prod_{i=1}^{\ell-k}(t-\xi_{i}) \prod_{j=\ell+1}^{n-1}(t-\xi_{j})
 + \beta_{n-1, k}\prod_{i=1}^{n-k-1}(t-\xi_{i}),\label{eq:eng2}
 \end{align}
 where $g_{n-k-1}$ is a polynomial of degree $n-k-1$.
We need $\beta_{kk},\ldots,\beta_{n-1,k}$ such that
$g_{n-k-1}$ is a linear combination
of the polynomials
\begin{equation}\label{eq:polw}
w_\ell(t):=\prod_{i=1}^{\ell-1}(t-\xi_{i}) \prod_{j=\ell+k}^{n-1}(t-\xi_{j}),\qquad
\ell=1,\ldots,n-k.
\end{equation}
The following lemma shows that such $\beta_{kk},\ldots,\beta_{n-1,k}$ do exist because
the above polynomials form a basis for the set of polynomials of
degree less than or equal to $n-k-1$.
\medskip

\begin{lem}\label{lem:nonsm} Polynomials $w_1,\ldots,w_{n-k}$ defined in
\eqref{eq:polw} form a basis for all polynomials
of degree less than or equal to $n-k-1$ for $k=1,\ldots,n-2$.
\end{lem}
\begin{pf}
We have $n-k$ polynomials $w_1,\ldots,w_{n-k}$ of degree $n-k-1$. If we look at
the matrix of values of these polynomials in points $\xi_1,\ldots,\xi_{n-k}$,
\[\left[\begin{matrix}w_1(\xi_1) & w_2(\xi_1) & \cdots & w_{n-k}(\xi_1)\cr
w_1(\xi_2) & w_2(\xi_2) & \cdots & w_{n-k}(\xi_2)\cr
\vdots &  \vdots & & \vdots \cr
w_1(\xi_{n-k}) & w_2(\xi_{n-k}) & \cdots & w_{n-k}(\xi_{n-k})\end{matrix}\right],
\]
we see that the matrix is lower triangular
with nonzero elements on the diagonal, because it follows
 from \eqref{eq:polw} that $w_\ell(\xi_j)=0$ for $j<\ell$ and $w_\ell(\xi_\ell)\ne 0$ for
 $\ell=1,\ldots,n-k$. The matrix
is therefore nonsingular and
the polynomials $w_1,\ldots,w_{n-k}$ satisfy the Haar condition. As a result $w_1,\ldots,w_{n-k}$
form a basis for all polynomials of degree less than or equal to $n-k-1$. \qed
\end{pf}

It is well known that two polynomials of degree $n-k-1$ or less are equal if they have
the same values at $n-k$ distinct points. We use this to write down a
system of linear equations for $\beta_{kk},\ldots,\beta_{n-1,k}$ from \eqref{eq:eng1} and \eqref{eq:eng2} by choosing
$\xi_1,\ldots,\xi_{n-k}$ as $n-k$ distinct points. This gives
\begin{equation}\label{eq:sysbeta}
\left[\begin{matrix}w_1(\xi_1) &  & & \cr
w_1(\xi_2) & w_2(\xi_2) &\cr
\vdots &  & \ddots & \cr
w_1(\xi_{n-k}) & w_2(\xi_{n-k}) & \cdots & w_{n-k}(\xi_{n-k})\end{matrix}\right]
\left[\begin{matrix} \beta_{kk} \cr \beta_{k+1,k} \cr \vdots \cr \beta_{n-1,k}\end{matrix}\right]
=
\left[\begin{matrix} g_{n-k-1}(\xi_1) \cr g_{n-k-1}(\xi_2) \cr \vdots \cr g_{n-k-1}(\xi_{n-k})\end{matrix}\right]
\end{equation}
and we know from Lemma \ref{lem:nonsm} that the system is nonsingular.
In addition to being lower triangular the system is also banded, and thus  can be solved  efficiently
with the forward substitution.

Once we set the coefficients $f_{kk},\ldots,f_{n-1,k}$
the updated polynomials $q_0^{(k)},\ldots,q_{n-1}^{(k)}$ give the
intermediate residual
\begin{align*}
s^{(k)}(x,y)&:=p(x,y)-p_{n0}x q_{n-1}^{(k)}(x,y)-\gamma_{n-2}q_{n-2}^{(k)}(x,y)-\cdots-
 \gamma_{n-k+1}q_{n-k+1}^{(k)}(x,y)\\
&\phantom{:}=s^{(k)}_{00}+s^{(k)}_{10}x+s^{(k)}_{01}y+\cdots+s^{(k)}_{n-k,0}s^{n-k}+
  \cdots+s^{(k)}_{0,n-k}y^{n-k}.
\end{align*}
We now insert $\gamma_{n-k}q_{n-k}^{(k)}$ to
zero the coefficient at $y^{n-k}$.
A simple computation shows that we have to choose
\begin{equation}\label{eq:gamma}
\gamma_{n-k}=(-1)^{n-k}\frac{s^{(k)}_{0,n-k}}{\xi_1\ldots \xi_{n-k}},
\end{equation}
which is well defined because $\xi_1,\ldots,\xi_{n-2}\ne 0$. The new
residual is
\begin{align*}
r^{(k)}(x,y)&:=p(x,y)-p_{n0}x q_{n-1}^{(k)}(x,y)-\gamma_{n-2}q_{n-2}^{(k)}(x,y)-\cdots-
 \gamma_{n-k}q_{n-k}^{(k)}(x,y)
\end{align*}
and compared to \eqref{eq:respk1} the degree of the residual is reduced by one.
We repeat the procedure and after $n-1$ steps we get the final residual
$r^{(n-1)}(x,y)=r^{(n-1)}_{00}+r^{(n-1)}_{10}x.$
The overall algorithm is presented in Algorithm \ref{alg:osijek}.
\bigskip

\begin{alg}\rm
Given a bivariate
polynomial of degree $n$
\[p(x,y)=p_{00}+p_{10}x+p_{01}y+\cdots + p_{n0}x^n+p_{n-1,1}x^{n-1}y +\cdots
+p_{1,n-1}x y^{n-1}\]
such that $p_{n0}\ne 0$, $p_{0n}=p_{0,n-1}=0$,
and all roots of $p_{n0}\xi^{n-1}+p_{n-1,1}\xi^{n-2} +\cdots + p_{1,n-1}=0$ are simple
and nonzero,
the output are $n\times n$ matrices $A$, $B$, and $C$ such that
$\det(A+xB+yC)=p(x,y)$.\label{alg:osijek}
\bigskip

\noindent\hbox{}\ \ 1. Compute the roots $\xi_1,\ldots,\xi_{n-1}$ of
$p_{n0}\xi^{n-1}+p_{n-1,1}\xi^{n-2} +\cdots + p_{1,n-1}=0$.\\[0.2em]
\hbox{}\phantom{1}2. $q_0(x,y)=1$\\[0.2em]
\hbox{}\phantom{1}3. for $j=1,\ldots,n-1$\\[0.2em]
\hbox{}\phantom{1}4. \qquad $q_j(x,y) = (x-\xi_jy)q_{j-1}(x,y)$\\[0.2em]
\hbox{}\phantom{1}5. $r(x,y)=p(x,y)-p_{n0} x q_{n-1}(x,y)$\\[0.5em]
\hbox{}\phantom{1}6. for $k=2,\ldots,n-1$\\[0.2em]
\hbox{}\phantom{1}7.\qquad $h(t)= (r_{n-k+1,0}t^{n-k} + r_{n-k,1}t^{n-k-1} + \cdots + r_{1,n-k} )/p_{n0}$\\[0.2em]
\hbox{}\phantom{1}8.\qquad $\alpha_{n-1,k} = r_{n-k+1,0}/p_{n0}-n+k+1$\\[0.2em]
\hbox{}\phantom{1}9.\qquad $g(t)=h(t)- t
\sum_{\ell=k}^{n-2}\prod_{i=1}^{\ell-k}(t-\xi_{i}) \prod_{j=\ell+1}^{n-1}(t-\xi_{j})
-\alpha_{n-1, k} t\prod_{i=1}^{n-k-1}(t-\xi_{i})$\\[0.2em]
\hbox{}10.\qquad for $m=1,\ldots,n-k$\\[0.2em]
\hbox{}11.\qquad\qquad for $\ell=1,\ldots,m$\\[0.2em]
\hbox{}12.\qquad\qquad\qquad $w_{m \ell}=\prod_{i=1}^{\ell-1}(\xi_{m}-\xi_{i}) \prod_{j=\ell+k}^{n-1}(\xi_{m}-\xi_{j})$\\
\hbox{}13.\qquad Solve
$\left[\begin{matrix}w_{11} &  &  & \cr
w_{21} & w_{22} &  & \cr
\vdots &  & \ddots & \cr
w_{n-k,1} & w_{n-k,2} & \cdots & w_{n-k,n-k}\end{matrix}\right]
\left[\begin{matrix} \beta_{kk} \cr \beta_{k+1,k} \cr \vdots \cr \beta_{n-1,k}\end{matrix}\right]
=
\left[\begin{matrix} g(\xi_1) \cr g(\xi_2) \cr \vdots \cr g(\xi_{n-k})\end{matrix}\right].
$\\[0.2em]
\hbox{}14.\qquad $q_k(x,y) = q_k(x,y)+ (x + \beta_{kk}y)q_{0}(x,y)$\\[0.2em]
\hbox{}15.\qquad for $j=k+1,\ldots,n-2$\\[0.2em]
\hbox{}16.\qquad\qquad $q_j(x,y) = (x-\xi_jy)q_{j-1}(x,y)+ \sum_{\ell=2}^k (x + \beta_{j\ell}y)q_{j-\ell}(x,y)$\\[0.2em]
\hbox{}17.\qquad $q_{n-1}(x,y) = (x-\xi_{n-1}y)q_{n-2}(x,y)+ \sum_{\ell=2}^k (\alpha_{n-1,\ell}x + \beta_{n-1,\ell}y)q_{n-1-\ell}(x,y)$\\[0.2em]
\hbox{}18.\qquad $s(x,y)=p(x,y)-p_{n0} x q_{n-1}(x,y) -\sum_{\ell=2}^{k-1}\gamma_{n-\ell}q_{n-\ell}(x,y)$\\[0.2em]
\hbox{}19.\qquad $\gamma_{n-k}=(-1)^{n-k}s_{0,n-k}/(\xi_1\ldots \xi_{n-k})$\\[0.2em]
\hbox{}20.\qquad $r(x,y)=s(x,y)-\gamma_{n-k}q_{n-k}(x,y)$\\[0.5em]
Return\\[0.2em]
$A=\left[\begin{matrix}r_{00} & \gamma_1 & \cdots & \gamma_{n-2} & 0 \cr
 & 1 \cr
 &  & \ddots\cr
 & &  & 1 \cr
 &  &  &  & 1
\end{matrix}\right],\qquad
B=\left[\begin{matrix} r_{10} & 0 & \cdots & 0 & p_{n0} \cr
-1 & 0\cr
\vdots &\ddots & \ddots \cr
-1 & \cdots & -1 & 0\cr
-\alpha_{n-1,n-1} & \cdots & -\alpha_{n-1,2} & -1 & 0
\end{matrix}\right],$\\[0.2em]
$C=\left[\begin{matrix}0 & & & & \cr
-\beta_{11} & 0 \cr
-\beta_{22} & -\beta_{21} & 0\cr
\vdots & & \ddots & \ddots \cr
-\beta_{n-1,n-1} & -\beta_{n-1,n-2} & \cdots & -\beta_{n-1,1} & 0
\end{matrix}\right].$
\end{alg}
\bigskip

\begin{rem}\label{rem:double} In Algorithm \ref{alg:osijek} we require that
 all roots of \eqref{eq:racxi} are simple and nonzero.  The nonzero condition is related to
expression \eqref{eq:gamma} for $\gamma_{n-k}$ that involves only
the roots $\xi_1,\ldots,\xi_{n-2}$. Thus it
is sufficient that the roots of \eqref{eq:racxi} are simple. If
one of the roots is zero we order them so that $\xi_{n-1}=0$
and then everything works fine.
\end{rem}

To further
clarify the algorithm, we give Example \ref{ex:po1l} with most of the details for the construction
of a determinantal representation of a bivariate polynomial of degree 5.
\bigskip

\begin{exm}\rm\label{ex:po1l}
We would like to linearize the bivariate polynomial
\begin{align*}p(x,y)&:=1 - x - 3 y + 3 x^2 - 7 x y - 6 y^2 + 10 x^3 + 9 x^2 y  - 14 x y^2 - 4 y^3\\
&\quad \ + 8 x^4 + 7 x^3 y - 8 x^2 y^2  - 4 x y^3 + 2 x^5 - 10x^3 y^2 + 8 x y^4
 \end{align*}
that satisfies the conditions $p_{50}\ne 0$ and $p_{05}=p_{04}=0$.

In the initial part before the main loop we  compute the roots $\xi_1,\ldots,\xi_4$ of the polynomial
\[p_{50}\xi^4+p_{41}\xi^3+p_{32}\xi^2+p_{23}\xi+p_{14}
 = 2\xi^4-10\xi^2+8.\]
The roots, which are all simple and nonzero, are
\[\xi_1=-2,\quad \xi_2=-1,\quad \xi_3=1,\quad \xi_4=2.\]
This gives the coefficients $f_{11},\ldots,f_{41}$ in the first subdiagonal
of \eqref{eq:ansatz}:
\[
f_{11}=x+2y,\quad f_{21}=x+y, \quad f_{31}=x-y, \quad f_{41}=x-2y.
\]
The coefficient at $q_4$ in the first row of \eqref{eq:ansatz} is $2x$.
The corresponding residual  is
\begin{align*}
r^{(1)}(x,y)
&=1 - x - 3 y + 3 x^2 - 7 x y - 6 y^2 + 10 x^3 + 9 x^2 y  - 14 x y^2 - 4 y^3\\
&\phantom{=}+ 8 x^4 + 7 x^3 y - 8 x^2 y^2  - 4 x y^3.
\end{align*}

Now we enter  the main loop.

\begin{itemize}
\item ($k=2$)
From the terms of $r^{(1)}$ of degree $4$ we define
\[s_4(x,y)=8 x^4 + 7 x^3 y - 8 x^2 y^2  - 4 x y^3\]
and divide $s_4(t,1)$ by $2t$ to obtain
$h_3(t)=4 t^3 + \frac{7}{2} t^2 - 4 t  - 4$.
Now we have to find the coefficients $\beta_{22},\beta_{32},\beta_{42}$, and $\alpha_{42}$ such that
\[h_3(t)=
(t+\beta_{22})(t-\xi_3)(t-\xi_4)+(t+\beta_{32})(t-\xi_1)(t-\xi_4)+(\alpha_{42}t+\beta_{32})(t-\xi_1)(t-\xi_2).\]
We get $\alpha_{42}$ from \eqref{eq:alpha} as
$\alpha_{42}=4-2=2$. This gives
\[g_2(t)=h_3(t)-t(t-\xi_3)(t-\xi_4)-t(t-\xi_1)(t-\xi_4)-\alpha_{42}t(t-\xi_1)(t-\xi_2)=
\frac{1}{2}t^2-6t-2.\]
 For $\beta_{22},\beta_{32}$, and $\beta_{42}$ we set the linear system
\[\left[\begin{matrix}
(\xi_1-\xi_3)(\xi_1-\xi_4) & 0 & 0 \cr
(\xi_2-\xi_3)(\xi_2-\xi_4) & (\xi_2-\xi_1)(\xi_2-\xi_4) & 0\cr
0 & (\xi_3-\xi_1)(\xi_3-\xi_4) & (\xi_3-\xi_1)(\xi_3-\xi_2)
\end{matrix}\right]\left[\begin{matrix}\beta_{22}\cr \beta_{32}\cr \beta_{42}\end{matrix}\right]=
\left[\begin{matrix}g_2(\xi_1) \cr g_2(\xi_2) \cr g_2(\xi_3)\end{matrix}\right].\]
When we insert the values we get
\[\left[\begin{matrix}
12 & 0 & 0 \cr
6 & -3 & 0\cr
0 & -3 & 6
\end{matrix}\right]\left[\begin{matrix}\beta_{22}\cr \beta_{32}\cr \beta_{42}\end{matrix}\right]=
\left[\begin{matrix}12 \\[0.2em] \frac{9}{2} \\[0.2em] -\frac{15}{2}\end{matrix}\right]\]
and the solution
is $\beta_{22}=1$, $\beta_{32}=\frac{1}{2}$, and $\beta_{42}=-1$,
which gives
\begin{align*}
s^{(2)}(x,y)
&=1 - x -  3 y + 3 x^2 - 7 x y -6 y^2 + 6 x^3 + 7 x^2 y  -
 12 x y^2 + 4 y^3.
\end{align*}
We compute $\gamma_3$ from \eqref{eq:gamma}
  to annihilate the term $y^3$ in $s^{(2)}$. We get
$\gamma_3 =  -4/(\xi_1\xi_2\xi_3)= 2$.
The residual at the end of step $k=2$ is
\begin{align*}
r^{(2)}(x,y)
&=1 - x -  3 y - x^2 - 12 x y - 6 y^2 + 4 x^3 + 3 x^2 y  -
 10 x y^2.
\end{align*}

\item ($k=3$) From $r^{(2)}$ we get
 $h_2(t)=2t^2+\frac{3}{2}t-5$.
Now we need $\beta_{33},\beta_{43}$, and $\alpha_{43}$ such that
\[h_2(t)=
(t+\beta_{33})(t-\xi_4)+(\alpha_{43}t+\beta_{43})(t-\xi_1).\]
We get $\alpha_{43}$ from \eqref{eq:alpha} as
$\alpha_{43}=2-1=1$. This gives
\[g_1(t)=h_2(t)-t(t-\xi_4)-\alpha_{43}t(t-\xi_1)=
\frac{3}{2}t-5.\]
 For $\beta_{33}$ and $\beta_{43}$ we set the linear system
\[\left[\begin{matrix}
\xi_1-\xi_4 & 0 \cr
\xi_2-\xi_4 & \xi_2-\xi_1
\end{matrix}\right]\left[\begin{matrix}\beta_{33}\cr \beta_{43}\end{matrix}\right]=
\left[\begin{matrix}g_1(\xi_1) \cr g_1(\xi_2)\end{matrix}\right].\]
When we insert the values we get
\[\left[\begin{matrix}
-4 & 0  \cr
-3 & 1
\end{matrix}\right]\left[\begin{matrix}\beta_{33}\cr \beta_{43}\end{matrix}\right]=
\left[\begin{matrix}-8 \\[0.2em] -\frac{13}{2} \end{matrix}\right]\]
and the solution is $\beta_{33}=2$ and $\beta_{43}=-\frac{1}{2}$.
This
gives
\[
s^{(3)}(x,y)=
1 -3 x -7 y - x^2 - 12 x y - 6 y^2.
\]
We compute $\gamma_2 =  -6/(\xi_1\xi_2)= -3$ from \eqref{eq:gamma}
to annihilate the term $y^2$ in $s^{(3)}$.
The residual after step $k=3$ is
\begin{align*}
r^{(3)}(x,y)
&=1 - 4 y +2 x^2 - 3 x y.
\end{align*}

\item ($k=4$) From $r^{(3)}$ we get
 $h_1(t)=t-\frac{3}{2}$.
Now we have to find coefficients $\beta_{44}$ and $\alpha_{44}$ such that
$h_1(t)=
\alpha_{44} t +\beta_{44}$. Clearly, the answer is
$\alpha_{44}=1$ and $\beta_{44}=-\frac{3}{2}.$
The new residual is
\[ s^{(4)}(x,y)
=1 -4y.\]
From \eqref{eq:gamma}
we compute
$\gamma_1 =  4/\xi_1= -2$.
The final residual after the main loop is
\begin{align*}
r^{(4)}(x,y)
&=1 + 2x.
\end{align*}
\end{itemize}

The final $5\times 5$ determinantal representation is
\[A+xB+yC=\left[\begin{matrix}
1+2x & -2 & 3 & 2 & 2x \\[0.2em]
-x-2y & 1 \\[0.2em]
-x-y & -x-y & 1 \\[0.2em]
-x-2y & -x-\frac{1}{2}y & -x + y & 1 \\[0.2em]
-x+\frac{3}{2}y & -x+\frac{1}{2}y & -2x+y & -x+2y & 1\end{matrix}\right].
\]
\end{exm}

Let us remark that the polynomial in Example \ref{ex:po1l} was constructed in such way that the
roots $\xi_1,\ldots,\xi_{n-1}$ are all real, which results in real matrices
in the determinantal representation. In a generic case, even if the
polynomial is real, the roots can be
complex and the representation has complex matrices.

\section{Determinantal representation for a non square-free polynomial}\label{sec:three}

If a bivariate polynomial $p$ is not square-free, then it
cannot be transformed into a polynomial that satisfies the conditions of Algorithm \ref{alg:osijek}. In such case we have several options.

First, for a polynomial of degree $n\le 5$ we can apply the algorithm in \cite{Anita} that returns a determinantal representation of order $n$ for a non square-free polynomial as well. For a polynomial of degree $n>5$ we can use one of the available symbolic or numerical tools, e.g., NAClab \cite{NAClab}, and factorize  $p$
into a product $p(x,y)=p_1(x,y)p_2(x,y)\ldots p_k(x,y)$, where $p_i$ is a square-free polynomial of degree $d_i$ for $i=1,\ldots,k$ and $d_1+\cdots+d_k=n$.  Now we apply
 Algorithm \ref{alg:osijek} to obtain $d_i\times d_i$ matrices $A_i,B_i$, and $C_i$ such that
$p_i(x,y)=\det(A_i+x B_i + y C_i)$ for $i=1,\ldots,k$. We arrange them in
block diagonal matrices $A=\diag(A_1,\ldots,A_k)$, $B=\diag(B_1,\ldots,B_k)$,
and $C=\diag(C_1,\ldots,C_k)$ and get $n\times n$ matrices such that
$\det(A+xB+yC)=p(x,y)$. So, combined with a square-free factorization, one can
find a determinantal representation of order $n$ for each bivariate polynomial
of degree $n$.

As the square-free factorization is more complex than Algorithm \ref{alg:osijek}, it takes most of the computational time in the above procedure. In our case, where we use representations to compute roots of a system of two bivariate polynomials, it is
more efficient for polynomials of small degree to use larger representations that can be constructed faster. For instance, for each polynomial of degree $n$ there exists a
uniform determinantal representation of order $2n-1$ \cite{Jan}. As an example,
a uniform representation of order $7$ for the polynomial \eqref{eq:pol} of degree $n=4$ is
\[
A+xB+yC=-
\begin{bmatrix}
-x     & \ph-1  &      &      &      &     &         \\
      &   -x  & \ph-1  &      &      &     &         \\
      &      &   -x  & 1  &      &     &         \\
p_{00}&p_{10}&p_{20}&p_{30}+p_{40} x&  -y &      \\
p_{01}&p_{11}&p_{21}+p_{31}x & & \ph-1 &  -y &   \\
p_{02}+ p_{03} y&p_{12}+c_{22}x& &     &      & \ph-1 &   -y     \\
 p_{13} x+p_{04} y& &     &      &      &     & \ph-1 \\
\end{bmatrix}.
\]
The construction of a uniform representation of order $2n-1$ is immediate, but we need more time in the second phase to numerically solve a larger associated two-parameter eigenvalue problem, which is singular in addition. Still, for polynomials of small degree this might altogether be faster than using a representation of order $n$ that requires a square-free factorization.

\section{Numerical examples}\label{sec:numex}

A new numerical approach for computing roots of systems of bivariate
polynomials by
constructing determinantal representations of the polynomials and solving the obtained two-parameter
eigenvalue problem
was proposed in \cite{BorMichiel}.
Using the representation of order $\frac{1}{6}n^2+{\cal O}(n)$ for a bivariate polynomial of degree $n$ (we refer to
this method as {\tt Lin2}) it was shown
that the approach is competitive for polynomials of
degree $9$ or less.
The method was compared to several
numerical methods for polynomial systems:
{\tt NSolve} in Mathematica~9~\cite{Wolfram},
BertiniLab~1.4~\cite{Bertini} running Bertini~1.5~\cite{BertiniExe},
NAClab~3.0~\cite{NAClab},
and {\tt PHCLab}~1.04~\cite{PHClab} running PHCpack~2.3.84,
which was the fastest of these methods.

Using the uniform representation of a smaller order $2n-1$ for
a bivariate polynomial of degree $n$ (we
refer to this method as {\tt MinUnif}),  computational times
were significantly improved in \cite{Jan} and the largest
degree such that the approach is competitive with
the existing numerical methods for polynomial systems was raised to 15.

In a numerical experiment we compare the new representation of order $n$
from Section \ref{sec:two}
(we refer to it as {\tt MinRep}) to {\tt Lin2}, {\tt MinUnif}, and {\tt PHCLab} on a similar set of random polynomials
as in \cite{Jan} and \cite{BorMichiel}. We take systems of
full bivariate polynomials of the same degree, whose coefficients are random real
numbers uniformly distributed on $[0,1]$ or random complex numbers, such that real
and imaginary parts are both uniformly distributed on $[0,1]$.

In order to overcome some difficulties that we noticed while testing  {\tt MinRep}, some heuristics were applied in the implementation. In practice,
even if the polynomial satisfies the initial conditions
of
Algorithm \ref{alg:osijek}, i.e, the roots of
\eqref{eq:racxi} are nonzero and simple, the obtained matrices
$A$, $B$, and $C$ can be such that the error between
$p(x,y)$ and $\det(A+xB+yC)$, caused by the numerical computation,
is too large. This happens for instance when the roots of
\eqref{eq:racxi} are ill-conditioned. Also, when the roots of \eqref{eq:racxi} are
close to each other, then linear systems \eqref{eq:sysbeta} can be
ill-conditioned.
A usual remedy for this is to
apply a random linear transformation \eqref{eq:trt}.

To get out of such troubles, we compute the determinantal
representation and check its quality by computing
\[\nu:=\max_{i=1,\ldots,k}\frac{|p(x_i,y_i)-\det(A+x_iB+y_iC)|}{
|p(x_i,y_i)|+\epsilon}\]
on a set of $k$ random points $(x_i,y_i)$, $i=1,\ldots,k$.
If $\nu\cdot \max(\|A\|_\infty,\|B\|_\infty, \|C\|_\infty)>\delta$ for
a given $\delta$ (in our experiments we use
$\epsilon=10^{-4}$, $k=200$, and $\delta=10^{-8})$, then we first compute a new representation of the polynomial, where we exchange the roles of
$x$ and $y$, and, if this does not help, then we apply a random change of variables \eqref{eq:trt}. We noticed
that this heuristics does not improve the situation for polynomials of degree 11 or more.
It seems that for a generic bivariate polynomial the largest safe degree,
when we can expect that the method works, is 10. This does not mean that the method cannot fail for polynomials
of smaller degree. Similar to other methods based on determinantal representations, this can happen for some systems that appear to be to
difficult for this approach and require computation in higher
precision.

\begin{table}[!htbp]
\begin{footnotesize}
\begin{center}
\caption{Average computational times in milliseconds
for {\tt MinRep}, {\tt MinUnif}, and {\tt PHCLab} for
random bivariate polynomial systems of degree $3$ to $10$. For  {\tt MinUnif}
separate results are included for real $(\RR)$ and complex polynomials $(\CC)$.
}\label{tbl:prim}
\begin{tabular}{r|rrrrr} \hline \rule{0pt}{2.3ex}%
$d$ & {\tt MinRep} & {\tt MinUnif} ($\RR$) & {\tt MinUnif} ($\CC$) & {\tt Lin2} & {\tt PHCLab}  \\
\hline \rule{0pt}{2.3ex}%
3 & 6 & 6 & 6 & 7 & 210 \\
4 & 8 & 9 & 11 & 11 & 247 \\
5 & 11 & 15 & 18 & 18 & 289 \\
6 & 15 & 25 & 32 & 32 & 344 \\
7 & 20 & 40 & 55 & 70 & 409 \\
8 & 29 & 70 & 98 & 191 & 499 \\
9 & 46 & 112 & 172 & 439 & 607 \\
10 & 64 & 184 & 301 & 1111 & 739 \\
\hline
\end{tabular}
\end{center}
\end{footnotesize}
\end{table}

The results in Table \ref{tbl:prim} show that the new
representation has a big potential as it is much faster than the previous ones. The results were obtained on a 64-bit Windows version of
MATLAB R2015b running on an Intel Core i5-6200U 2.30 GHz processor with 8 GB
of RAM.
For each $n$ we apply the methods to the same set of
50 real and 50 complex random polynomial systems of degree $n$
and measure the average time. The accuracy of all methods on this
random set of polynomials is comparable.
For {\tt MinUnif}, where determinantal representations have real
matrices for real polynomials, we report separate results
for polynomials with real and complex coefficients.
Although {\tt MinRep} is still the fastest method
for $11\le n\le 15$, the computed roots are in most cases
 useless and, based on the results
from \cite{Jan}, we rather suggest that {\tt MinUnif} is applied
to polynomials of such degree. As in {\tt MinUnif} no computation is involved
in the construction of determinantal representations,
this phase is not affected by numerical errors. On the other hand, the matrices
in {\tt MinUnif} are approximately double size compared to
the matrices in {\tt MinRep} and to extract the final solution
one has to compute the finite eigenvalues of a pair of singular pencils,
which is more delicate than solving a regular pencil with the QZ algorithm
in {\tt MinRep}.
A future research might give more insight into the cases when {\tt MinRep}
does not perform so well and further improve the method.

We performed a limited number of  tests comparing {\tt MinUnif}
to the square-free factorization approach from Section \ref{sec:three}. We observed that the square-free factorization approach is faster for polynomials of degree
$n\ge 8$. Although
{\tt MinUnif} is faster for polynomials of degree $n=6$ and $n=7$, it is also less
accurate because of the multiple roots. We therefore recommend to use
the square-free factorization for all non square-free polynomials of degree $n\ge 5$.

\section{Conclusions}\label{sec:conc}

We presented a numerical construction for the determinantal representation
of a square-free bivariate polynomial of degree $n$ with matrices of order $n$.
The computation requires only  routines for roots of
univariate polynomials of degree $n$ and solutions of linear systems of order
less than $n$, both should be available in any numerical package.
When combined with a square-free factorization, we believe that this is
the first construction that reaches the theoretical lower
order from Dixon's theorem for
all bivariate polynomials. We do not get symmetric matrices,
but this is not really important in our application.

The new representation can be used to numerically compute the roots
of a system of bivariate polynomials. Compared to the previous results,
the new representation speeds up the computation considerably.

\bigskip\noindent{\bf Acknowledgment } The author wishes to thank
Anita Buckley (University of Ljubljana) and Michiel~E. Hochstenbach (TU Eindhoven)
for helpful discussions.
Part of the research was performed while
the  author was visiting the Department of Mathematics at University of Osijek.
The author wishes to thank Ninoslav Truhar (University of Osijek) for the hospitality.
The author would also like to thank the referees for careful reading of the paper
and several helpful suggestions.
\medskip

The author has been partially supported by the
Slovenian Research Agency (research core funding P1-0294).


\begin{thebibliography}{99}

\bibitem{Atkinson} {F.V.~Atkinson},
  {Multiparameter Eigenvalue Problems}, Academic Press, New York, 1972.

\bibitem{BertiniExe} {D.J.~Bates, J.H.~Husenstein, A.J.~Sommese, C.W.~Wampler},
  {Bertini: Software for Numerical Algebraic Geometry}, available at {\tt bertini.nd.edu};
  doi: {\tt dx.doi.org/10.7274/R0H41PB5}.
\bibitem{Jan} {A.~Boralevi, J.~van~Doornmalen, J.~Draisma,
M.E.~Hochstenbach, B.~Plestenjak}, Uniform
determinantal representations, arXiv:1607.04873 (2016), to appear in SIAM
Journal on Applied Algebra and Geometry.
\bibitem{Anita} {A.~Buckley, B.~Plestenjak}, Simple determinantal representations of up to quintic bivariate polynomials, arXiv.1609.00498 (2016).
\bibitem{Dixon} {A.~Dixon}. {Note on the reduction of a ternary quartic to a symmetrical determinant.}
  Proc.~Camb.~Phil.~Soc.~11 (1902) 350--351.

\bibitem{PHClab} {Y.~Guan, J.~Verschelde}, PHClab: A MATLAB/Octave interface to PHCpack.
In: M.~Stillman, J.~Verschelde, N.~Takayama (eds), {Software for Algebraic Geometry},
volume 148 of The IMA Volumes in Mathematics and its Applications, Springer,
New York, (2008) 15--32.


\bibitem{Khazanov} {V.B.~Khazanov},
 To solving spectral problems for multiparameter polynomial matrices,
 J.~Math.~Sci.~141 (2007) 1690--1700.

\bibitem{Matlab} {The MathWorks, Inc.}, Matlab, Natick, Massachusetts, United States.

\bibitem{MuhicPlestenjakLAA} {A.~Muhi\v{c}, B.~Plestenjak},
  {On the quadratic two-parameter eigenvalue problem and its linearization},
  Linear Algebra Appl.~432 (2010) 2529--2542.

\bibitem{Netzer} {T. Netzer, A.~Thom},
  {Polynomials with and without determinantal representations},
  Linear Algebra Appl.~437 (2012) 1579--1595.

\bibitem{Bertini} {A. Newell}, BertiniLab: toolbox for solving polynomial systems,
MATLAB Central File Exchange,
\url{www.mathworks.com/matlabcentral/fileexchange/48536-bertinilab}

\bibitem{Plaumann} {D.~Plaumann, R.~Sinn, D.E.~Speyer, C.~Vinzant},
 {Computing Hermitian determinantal representations of hyperbolic curves},
 Internat. J. Algebra Comput. 25 (2015) 1327--1336.

\bibitem{PlaumannHVcurves} {D.~Plaumann, B.~Sturmfels, C.~Vinzant}, {Computing Linear Matrix Representations of Helton-Vinnikov Curves}, Operator Theory: Advances and Applications 222 (2012) 259--277.

\bibitem{BorBR} {B.~Plestenjak},
  BiRoots, MATLAB Central File Exchange,\\
  \url{http://www.mathworks.com/matlabcentral/fileexchange/54159-biroots}.

\bibitem{BorMC1} { B.~Plestenjak},
  MultiParEig: toolbox for multiparameter eigenvalue problems, MATLAB Central File Exchange,
  \url{www.mathworks.com/matlabcentral/fileexchange/47844-multipareig}.


\bibitem{BorMichiel} {B.~Plestenjak, M.E.~Hochstenbach},
  Roots of bivariate polynomial systems via determinantal representations,
  SIAM J. Sci. Comput. 38 (2016) A765--A788.

\bibitem{Quarez} {R.~Quarez},
{Symmetric determinantal representation of polynomials},
Linear Algebra Appl.~436 (2012) 3642--3660.

\bibitem{Smith} {K.E.~Smith, L.~Kahanp\"a\"a, P.~Kek\"al\"ainen, W.~Traves}, {An
Invitation to Algebraic Geometry}, Springer, New York, 2000.

\bibitem{phc} {J.~Verschelde}, {Algorithm 795: PHCpack: a general-purpose solver for polynomial systems
by homotopy continuation}, ACM Trans.~Math.~Softw., 25 (1999) 251--276.

\bibitem{VinnikovSurvey} {V.~Vinnikov}, {LMI representations of convex
semialgebraic sets and determinantal representations of
algebraic hypersurfaces: past, present, and future},
in H.~Dym, M.C.~de Oliveira, M.~Putinar (eds.),
{Mathematical Methods in Systems, Optimization, and Control:
Festschrift in Honor of J.~William Helton.}
Operator Theory: Advances and Applications 222, Birkh\"auser (2012)
325--349.

\bibitem{Wolfram} {Wolfram Research, Inc.}, Mathematica, Version 9.0,
Champaign, Illinois, 2012.

 \bibitem{NAClab} {Z.~Zeng,  T.-Y.~Li}, NAClab: a Matlab toolbox for numerical
 algebraic computation,
 ACM Commun. Comput. Algebra~47 (2013) 170--173.

\end{thebibliography}
\end{document}